\documentclass[12 pt]{amsart}

\usepackage{amsthm}
\usepackage{amsmath}
\usepackage{amssymb}
\usepackage{latexsym}
\usepackage{url}
\usepackage{graphicx}
\usepackage{bmpsize}
\usepackage{mathtools}
\usepackage[margin=3cm]{geometry}

\author{Felix Goldberg}
\address{Caesarea-Rothschild Institute, University of Haifa, Haifa, Israel}
\email{felix.goldberg@gmail.com}

\title{Graph energy estimates via the Chebyshev functional}
\date{September 2, 2014}

\newtheorem{thm}{Theorem}[section]
\newtheorem{cor}[thm]{Corollary}

\newtheorem{conj}[thm]{Conjecture}

\newtheorem{rmrk}[thm]{Remark}

\DeclareMathOperator{\Tr}{Tr}

\begin{document}

\begin{abstract}
Let $G$ be a graph with $n$ vertices and $m$ edges. The energy $E$ of the graph $G$ is defined as the sum of the moduli of the adjacency eigenvalues $\lambda_{1} \geq \lambda_{2} \geq \ldots \geq \lambda_{n}$ of $G$:
$$
E=\sum_{i=1}^{n}{|\lambda_{i}|}.
$$
We obtain new lower bounds on the energy of a graph, which in various cases improve upon known results. For example, a particularly simple and appealing corollary of our results is:
$$
E \geq \frac{2m}{\lambda_{1}}.
$$
This implies a result obtained by Gutman \emph{et al.} for regular graphs and is better for triangle-free graphs than a result of Caporossi \emph{et al.}.
\end{abstract}



\thanks{{This research was supported by the Israel Science Foundation (grant number 862/10.)}}

\maketitle

\section{Introduction}
Let $G$ be a graph with adjacency matrix $A$ whose eigenvalues  are $\lambda_{1} \geq \lambda_{2} \geq \ldots \geq \lambda_{n}$. The energy of $G$ is then defined as:
$$
E=\sum_{i=1}^{n}{|\lambda_{i}|}.
$$
This concept has been first introduced and intensively studied in the context of mathematical chemistry but in the last 15 years it has garnered a lot of attention from graph theorists as well. For overviews of the subject we refer the reader to the recent book \cite{Energy} by Li, Shi, and Gutman and to the earlier surveys by Gutman \cite{Gutman_survey} and Brualdi \cite{Brualdi_energy}.

Our aim in this paper is to contribute a new lower bound for the energy, obtained by revisiting the original approach of one of the early pioneers, McClelland and bringing a discrete variant of the well-known Gr\"{u}ss integral inequality to bear on it.

\subsection{Notation, terminology, and some standard facts}
Throughout the paper we shall assume that the graph $G$ has $n$ vertices and $m$ edges. We shall denote by $t$ the smallest absolute value of an eigenvalue of $G$, that is: $t=\min\{|\lambda_{i}|\}$. The graph will be called \emph{singular} if $t=0$ and \emph{non-singular} otherwise.

By the Perron-Frobenius theorem it is known that $|\lambda_{1}| \geq |\lambda_{i}|$ for any $i$. We also have the following well-known fact:
$$
\sum_{i=1}^{n}{\lambda_{i}^{2}}=\Tr{A^{2}}=2m.
$$

\section{Some known results}\label{sec:sur}
In this brief section we do not attempt to provide an exhaustive survey, but rather to indicate the main lower bounds on energy that are present in the literature, so that the reader can compare them with the new result we shall derive. 

In the halcyon days of graph energy McClelland \cite{McCle71} obtained the following bounds:
\begin{thm}\cite{McCle71}\label{thm:mcc}
$$\sqrt{2m+n(n-1)|A|^{\frac{2}{n}}} \leq E \leq \sqrt{2mn}.$$
\end{thm}

A different lower bound has been given by Caporossi \emph{et al.} \cite{Caporossi_etal99}:
\begin{thm}\cite{Caporossi_etal99}\label{thm:caporossi}
$$E \geq 2\sqrt{m}.$$
\end{thm}

Clearly, for singular graphs Theorem \ref{thm:caporossi} is better than the lower bound of Theorem \ref{thm:mcc}.


\section{Statement of the new results and some discussion}\label{sec:main}
Our first main result is:
\begin{thm}\label{thm:main}
\begin{equation}\label{eq:main}
E \geq \frac{2m+n\lambda_{1}t}{\lambda_{1}+t}.
\end{equation}
\end{thm}

Since the right-hand side of \eqref{eq:main} is a non-decreasing function of $t$, we can deduce:
\begin{cor}\label{cor:nice}
$$E \geq \frac{2m}{\lambda_{1}}.$$
\end{cor}

\begin{rmrk}
Equality is attained in Corollary \ref{cor:nice} for the complete bipartite graphs $K_{p,q}$ with $E=2\sqrt{pq}$ and $m=2pq,\lambda_{1}=\sqrt{pq}$. I am grateful to Dr. Clive Elphick for this observation. It would be an interesting problem to try to find other graphs - if there are any - which attain equality.
\end{rmrk}

Since a $d$-regular graph has $2m=nd$ and $\lambda_{1}=d$, Corollary \ref{cor:nice} implies the following result by Gutman \emph{et al.} \cite{Gutman_etal07}:
\begin{thm}\cite{Gutman_etal07}
Let $G$ be an $r$-regular graph, $r>0$. Then $$E \geq n.$$
\end{thm}

Another class of graphs for which Corollary \ref{cor:nice} improves upon known results is that of triangle-free graphs. It is known \cite{Nosal} that for them $\lambda_{1} \leq \sqrt{m}$  and therefore Corollary \ref{cor:nice} is better than Theorem \ref{thm:caporossi} over this class.

It is also possible to deduce, using the arithmetic-geometric means inequality, another consequence of Theorem \ref{thm:main} which is not very strong but forms a nice counterpart to Theorem \ref{thm:mcc}:
\begin{cor}
$$
E \geq \sqrt{2mn} \cdot \sqrt{\frac{4\lambda_{1}t}{(\lambda_{1}+t)^{2}}}.
$$
\end{cor}

As we shall see, our approach will enable us to give an even stronger bound than \eqref{eq:main}. To state it we need to single out the smallest \emph{non-zero} eigenvalue of the graph:
$$
t_{nz}=\min\{|\lambda_{i}| \Big{\vert} \lambda_{i} \neq 0\}.
$$

\begin{thm}\label{thm:better}
Let $r=rank(A)$. Then
$$
E \geq \frac{2m+r\lambda_{1}t_{nz}}{\lambda_{1}+t_{nz}}.
$$
\end{thm}

We defer the proofs of Theorems \ref{thm:main} and \ref{thm:better} till Section \ref{sec:proof}, while in the next section we set up the engine of the proof.

\section{The Gr\"{u}ss inequality}
Chebyshev's classic inequality says that if $f,g:[a,b] \rightarrow \mathbb{R}$ are integrable functions, either both increasing or both decreasing, then
\begin{equation}\label{eq:cheb}
\int_{a}^{b}{f(x)g(x)dx} \geq \frac{1}{b-a}\int_{a}^{b}{f(x)dx}\int_{a}^{b}{g(x)dx}.
\end{equation}
This elegant inequality has been generalized and extended in many ways. We refer to \cite[Chapters IX--X]{MPF_book} for a survey of some of these developments.

In 1935 Gr\"{u}ss proved the following result:
\begin{thm}\cite{Gru35}\label{thm:gruss}
Let $f,g:[a,b] \rightarrow \mathbb{R}$ be integrable functions such that 
\begin{equation}\label{eq:gruss_cond}
\phi \leq f(x) \leq \Phi \textit{ and } \gamma \leq g(x) \leq \Gamma \textit{ for all } x \in [a,b].
\end{equation}
Then
$$
\Big|\frac{1}{b-a}\int_{a}^{b}{f(x)g(x)dx} - \frac{1}{b-a}\int_{a}^{b}{f(x)dx}\frac{1}{b-a}\int_{a}^{b}{g(x)dx}\Big| \leq \frac{1}{4}(\Phi-\phi)(\Gamma-\gamma).
$$

\end{thm}


Now let us state the abstract formulation due to Dragomir \cite{Dra99}, as it permits an easy derivation of the discrete variant we need: Let $(X,(\cdot))$ be a real inner product space and let $e \in X, ||e||=1$. The \emph{Chebyshev functional} on $X$ is defined as:
$$
\forall x,y \in X \quad T(x,y)=\langle x,y \rangle -\langle x,e \rangle \langle y,e \rangle.
$$
Taking the product with $e$ is the operation of taking a "mean". For any $z \in X$ we denote $A(z)=\langle x,e \rangle$.

\begin{thm}\cite{Dra99}\label{thm:drag1}
Let $x,y \in X$ be such that there exist $\phi,\gamma,\Phi,\Gamma \in \mathbb{R}$ so that the conditions
\begin{equation}\label{eq:gruss_abs}
\langle \Phi e-x,x-\phi e \rangle \geq 0 \textit{ and } \langle \Gamma e-y,y-\gamma e \rangle \geq 0 
\end{equation}
hold. Then 
$$
|T(x,y)| \leq \frac{1}{4}|\Phi-\phi||\Gamma-\gamma|.
$$
\end{thm}
Note that \eqref{eq:gruss_abs} reduces to \eqref{eq:gruss_cond} for the inner product $\langle f,g \rangle=\frac{1}{b-a}\int_{a}^{b}f(x)g(x)dx$.

We shall need a stronger version of Theorem \ref{thm:drag1}, which is implicit in Dragomir's proof:

\begin{thm}\label{thm:drag2}
Under the assumptions of Theorem \ref{thm:drag1},
$$
|T(x,y)| \leq \sqrt{|\Phi-A(x)||A(x)-\phi||\Gamma-A(y)||A(y)-\gamma|}.
$$
\end{thm}

Now let us equip $X=\mathbb{R}^{n}$ with the inner product $\langle x,y\rangle=\frac{1}{n}\sum_{i=1}^{n}{x_{i}y_{i}}$ and record the following consequence of Theorem \ref{thm:drag2}:
\begin{thm}\label{thm:drag3}
Let $x,y \in \mathbb{R}^{n}$ and let $A(x)=\frac{1}{n}\sum_{i=1}^{n}{x_{i}},A(y)=\frac{1}{n}\sum_{i=1}^{n}{y_{i}}$.
If $\phi \leq x_{i} \leq \Phi$ and $\gamma \leq y_{i} \leq \Gamma$, then
$$
\Big|\frac{1}{n}\sum_{i=1}^{n}{x_{i}y_{i}}-\frac{1}{n^{2}}\sum_{i=1}^{n}{x_{i}}\sum_{i=1}^{n}{y_{i}}\Big| \leq \sqrt{(\Phi-A(x))(A(x)-\phi)(\Gamma-A(y))(A(y)-\gamma)}.
$$
\end{thm}


\section{Proofs for Section \ref{sec:main}}\label{sec:proof}
The following observation goes back to McClelland \cite{McCle71}:
$$
E^{2}=\sum_{i=1}^{n}{|\lambda_{i}|^{2}}+\sum_{i \neq j}{|\lambda_{i}||\lambda_{j}|}=2m+\sum_{i \neq j}{|\lambda_{i}||\lambda_{j}|}.
$$
Let us denote $P=\sum_{i \neq j}{|\lambda_{i}||\lambda_{j}|}$. It is clear that estimating $E$ is equivalent to estimating $P$. In \cite{Caporossi_etal99} the bound $P \geq 2m$ was observed, leading to the claim of Theorem \ref{thm:caporossi}. We show here a different approach to bounding $P$, based on representing it as an inner product.

Let $x_{i}=|\lambda_{i}|$ and $y_{i}=E-|\lambda_{i}|$ for $1 \leq i \leq n$. Then it is easy to see that
$$
P=\sum_{i=1}^{n}{x_{i}y_{i}}.
$$
Observe that $$t \leq x_{i} \leq \lambda_{1} \textit{ and } E-\lambda_{1} \leq y_{i} \leq E-t.$$ Also, $$A(x)=\frac{E}{n} \textit{ and } A(y)=\frac{(n-1)E}{n}.$$

Now we apply Theorem \ref{thm:drag3} to obtain:
$$
\Big|\frac{P}{n}-\frac{(n-1)E^{2}}{n^{2}}\Big| \leq \sqrt{(\lambda_{1}-\frac{E}{n})(\frac{E}{n}-t)(\frac{E}{n}-t)(\lambda_{1}-\frac{E}{n})}=(\lambda_{1}-\frac{E}{n})(\frac{E}{n}-t).
$$
Therefore:
$$
P \geq n\Big(\frac{(n-1)E^{2}}{n^{2}}-(\lambda_{1}-\frac{E}{n})(\frac{E}{n}-t)\Big)=E^{2}+n\lambda_{1}t-(\lambda_{1}+t)E.
$$

This implies:
$$
E^{2}=2m+P \geq 2m+E^{2}+n\lambda_{1}t-(\lambda_{1}+t)E
$$
which immediately leads to \eqref{eq:main} upon trivial re-arrangement. This concludes the proof of Theorem \ref{thm:main}. \qed

\medskip
	
To prove Theorem \ref{thm:better} it is only necessary to observe that zero eigenvalues of $A$ correspond to zero entries in the vector $x$. Delete them and the corresponding entries from $y$ to obtain shorter vectors $x^{'},y^{'} \in \mathbb{R}^{r}$ which satisfy:
$$
P=\sum_{i=1}^{n}{x^{'}_{i}y^{'}_{i}},
$$
$$t_{nz} \leq x_{i} \leq \lambda_{1} \quad \quad  E-\lambda_{1} \leq y_{i} \leq E-t_{nz}$$ and $$A(x)=\frac{E}{r}  \quad \quad  A(y)=\frac{(r-1)E}{r}.$$	

Therefore the arguments given before work the same way, with $t$ replaced by $t_{nz}$ and $n$ replaced by $r$. \qed


\section{Two conjectures by Elphick}
Dr. Clive Elphick has communicated to me two very interesting conjectures engendered by the results reported here. To state them let us introduce two measures of the irregularity of a graph, studied in the paper \cite{ElpWoc13} by Elphick and Wocjan. Let $d_{1},d_{2},\ldots,d_{n}$ be the vertex degress of $G$ and let $d$ be the average degree. Define:
$$
\epsilon=\frac{n\sum_{i \sim j}{\sqrt{d_{i}d_{j}}}}{2m^{2}}
$$
and
$$
\beta=\frac{\lambda_{1}}{d}=\frac{\lambda_{1}n}{2m}.
$$

\smallskip
It is known that $\beta \geq \epsilon \geq 1$ (cf. \cite[p. 53]{ElpWoc13}). We can now state the conjectures.

\begin{conj}
Let $G$ be a connected graph. Then 
$$
E \geq \frac{n}{\epsilon}.
$$
\end{conj}
Since $\frac{2m}{\lambda_{1}}=\frac{n}{\beta}$, this would be an improvement of Corollary \ref{cor:nice}.

\begin{conj}
Let $G$ be a connected graph. Then 
$$
E \leq \frac{2m}{\sqrt{\lambda_{1}}}.
$$
\end{conj}
Since $\frac{2m}{\sqrt{\lambda_{1}}} \leq \sqrt{\frac{2mn}{\beta}} \leq \sqrt{2mn}$ this would an improvement over the upper bound of Theorem \ref{thm:mcc}.

Both conjectures have been verified for all $11117$ graphs on eight vertices using a computer. The connectedness assumption is essential for both conjectures.

\bibliographystyle{abbrv}
\bibliography{nuim}
\end{document}